\documentclass[12pt,english,reqno]{smfart}
\usepackage[letterpaper,left=1in,right=1in,top=1.5in,bottom=1.5in]{geometry}

\usepackage[T1]{fontenc}
\usepackage[english,francais]{babel}

\usepackage{amssymb,url,xspace,smfthm,float,subfig}

\usepackage{tikz-cd}

\usepackage[export]{adjustbox}

\newcommand{\SMF}{Soci\'et\'e Ma\-th\'e\-Ma\-ti\-que de France}
\newcommand{\BibTeX}{{\scshape Bib}\kern-.08em\TeX}
\newcommand{\T}{\S\kern .15em\relax }
\newcommand{\AMS}{$\mathcal{A}$\kern-.1667em\lower.5ex\hbox
        {$\mathcal{M}$}\kern-.125em$\mathcal{S}$}

\usepackage[hidelinks]{hyperref}
\usepackage[
]{cleveref}

\newcommand*{\TitleFont}{%
      \usefont{\encodingdefault}{\rmdefault}{}{n}%
      \fontsize{9}{20}%
      \selectfont}

\title[Knots, Primes and the adele class space]{Knots, Primes and the adele class space\\\TitleFont{}}\date {}
\author{Alain Connes}
\address{Coll\`ege de France\\3 Rue d'Ulm\\75005 Paris, France}
\address{I.H.E.S., France}
\email{\url{alain@connes.org}}
\author{Caterina Consani}
\address{Department of Mathematics\\Johns Hopkins University\\Baltimore, MD 21218, USA}
\email{\url{cconsan1@jhu.edu}}

\subjclass{\href{http://www.ams.org/mathscinet/msc/msc2020.html?t=11R37&btn=Current}{11R37},
\href{http://www.ams.org/mathscinet/msc/msc2020.html?t=11M06&btn=Current}{11M06},
\href{http://www.ams.org/mathscinet/msc/msc2020.html?t=11Mxx&btn=Current}{11M55},
\href{http://www.ams.org/mathscinet/msc/msc2020.html?t=14A15&btn=Current}{14A15},
\href{http://www.ams.org/mathscinet/msc/msc2020.html?t=14F20&btn=Current}{14F20},
\href{http://www.ams.org/mathscinet/msc/msc2020.html?t=57K10&btn=Current}{57K10}}

\usepackage{amsmath}
\usepackage{amscd}
\usepackage{amsbsy}
\usepackage{amssymb}
\usepackage{verbatim}
\usepackage{eufrak}
\usepackage{eucal}
\usepackage{microtype}
\usepackage{mathrsfs}
\usepackage{amsthm}
\usepackage{stmaryrd}
\usepackage[all,cmtip]{xy}
\usepackage{tikz}
\usetikzlibrary{matrix,arrows}
\usepackage[abbrev,lite,alphabetic]{amsrefs}
\usepackage{dsfont}

\usepackage{color}
\definecolor{todo}{rgb}{1,0,0}
\definecolor{conditional}{rgb}{0,1,0}
\definecolor{e-mail}{rgb}{0,.40,.80}
\definecolor{reference}{rgb}{.20,.60,.22}
\definecolor{mrnumber}{rgb}{.80,.40,0} 
\definecolor{citation}{rgb}{0,.40,.80} 

\def\qqq{\,,\,~\forall}

\def\hatz{{\hat\Z^\times}}
\def\Spec{{\rm Spec\,}}

\newcommand{\GL}{{\rm GL}}

\def\cA{\mathcal{A}}
\def\cO{{\mathcal O}}

\def\A{{\mathbb A}}
\def\C{{\mathbb C}}
\def\F{{\mathbb F}}
\def\Q{{\mathbb Q}}
\def\R{\mathbb{R}}
\def\Z{\mathbb{Z}}

\theoremstyle{plain}
\newtheorem{thm}{Theorem}[section]

\newtheorem{rem}[thm]{Remark}

\begin{document}
\def\smfbyname{}

\begin{abstract}
We show that the  scaling site $X_\Q$ and its periodic orbits $C_p$ of length $\log p$  offer a  geometric framework for the well-known analogy between primes and knots. The role of the maximal abelian cover of $X_\Q$ is played by the  quotient map $\pi:X_\Q^{ab}\to X_\Q$ from the adele class space $X_\Q^{ab}:=\Q^\times \backslash \A_\Q$ to   $X_\Q=X_\Q^{ab}/\hatz$.	 The inverse image $\pi^{-1}(C_p)\subset X_\Q^{ab}$ of the periodic orbit $C_p$  is canonically isomorphic to the mapping torus  of the multiplication by the Frobenius at $p$  in the abelianized \'etale fundamental group $\pi_1^{e t}(\Spec \, \Z_{(p)})^{ab}$ of the spectrum of the local ring $\Z_{(p)}$, thus exhibiting the linking of $p$ with all other primes.  In the same way as the Grothendieck theory of the  \'etale fundamental group of schemes is an extension of Galois theory to schemes, the adele class space gives, as a covering of the scaling site, the corresponding extension of the class field isomorphism  for $\Q$ to schemes related to $\Spec\,\Z$.
\end{abstract}

\maketitle

\tableofcontents

\section{Introduction}

The adele class space of \cite{Co-zeta} delivers both a geometric interpretation of the explicit Riemann-Weil formulas as a trace formula, as well as a spectral realization of the zeros of the Riemann zeta function and of $L$-functions with Gr\"ossencharakter. There is a well known analogy, due to Barry Mazur\footnote{on a suggestion of David Mumford},  between knots and primes \cite{Mazur0,Mazur,Mazur1,Mazur2,Morishita,Morishita1}. 
We show that the periodic orbits $C_p$ of length $\log p$ in the Riemann sector $X_\Q:=\Q^\times \backslash \A_\Q/\hatz$ of the adele class space provide a geometric realization of this conjectural relation. 
By construction $X_\Q$ is obtained as the quotient of $\Q^\times \backslash \A_\Q$ by the action of the compact subgroup $\hatz$ of idele classes, acting by multiplication. The space $X_\Q$ admits a simple description in terms of rank one groups \cite{CCas}. Namely, it is the union of the space of such groups, up to isomorphism, and the space of rank one subgroups of $\R$. This latter space is a prototype of  a noncommutative space.\newline
The action of the scaling group $\R_+^*$ on $X_\Q$ gives a Hasse-Weil interpretation of the Riemann zeta function \footnote{completed with its archimedean local factor}  \cite{CC1,ccmonoids} and coincides with the action of the Frobenius automorphisms on the points of the arithmetic site over the tropical semifield $\R_+^{max}$ \cite{CCas}. 
To each prime $p$ corresponds the  periodic orbit $C_p\subset X_\Q$ of length $\log p$ for the action of $\R_+^*$:  $C_p$ consists precisely of  rank one subgroups of $\R$ which are isomorphic to the additive group of the ring $\Z[1/p]$.  We refer to \cite{CCscal1} for the structure of $C_p$ inherited from the scaling site :  $C_p=\R_+^*/p^\Z$ appears as an  elliptic curve in characteristic one,  similar to the Jacobi elliptic curve $\C^*/q^\Z$. The Riemann-Roch formula for $C_p$ involves real valued dimensions. \newline
The role of  the maximal abelian cover $X_\Q^{ab}$ of $X_\Q$ is played by the adele class space $\Q^\times \backslash \A_\Q$ together with the  canonical quotient map $\pi$:
\begin{equation}\label{mappi}
\pi:X_\Q^{ab}:= \Q^\times \backslash \A_\Q\to X_\Q=\Q^\times \backslash \A_\Q/\hatz.
\end{equation}
 Let  $\Z_{(p)}$ be the  ring $\Z$ localized at $p$, $r:\Z_{(p)}\to \F_p$  the residue morphism, and $r^*$ the induced embedding of schemes
\begin{equation}\label{residue}
r^*:\Spec \, \F_p\hookrightarrow \Spec \, \Z_{(p)}.
\end{equation}
Our first result involves the map  $\pi_1^{e t}(r^*)$ at the level of  the \'etale fundamental groups $\pi_1^{e t}$ of schemes. One lets $F r o b_{p}$ be the canonical generator of $\pi_1^{e t}(\Spec\left(\F_p\right))\simeq{\rm Gal}_{}(\bar \F_p/\F_p)$ given by the Frobenius.\newline Theorem \ref{mainglobalintro} shows that the monodromy obtained by lifting the periodic orbit  $C_p$ in the maximal abelian cover $X_\Q^{ab}$ coincides with the multiplication by $r^*\left\{F r o b_{p}\right\}$ in the  abelianized \'etale fundamental group $\pi_1^{e t}(\Spec \, \Z_{(p)})$.
\begin{thm}\label{mainglobalintro} Let $p$ be a prime and $\left\{F r o b_{p}\right\}\in \pi_1^{e t}(\Spec\left(\F_p\right))$ be the canonical generator. The inverse image $\pi^{-1}(C_p)\subset X_\Q^{ab}$ of the periodic orbit $C_p$  is canonically isomorphic to the mapping torus\footnote{in the ordinary topological sense} of the multiplication by $r^*\left\{F r o b_{p}\right\}$  in the abelianized \'etale fundamental group $\pi_1^{e t}(\Spec \, \Z_{(p)})^{ab}$. The canonical isomorphism is equivariant for the action of the idele class group.\end{thm}
 \begin{figure}[H]
\begin{center}
\includegraphics[scale=0.35]{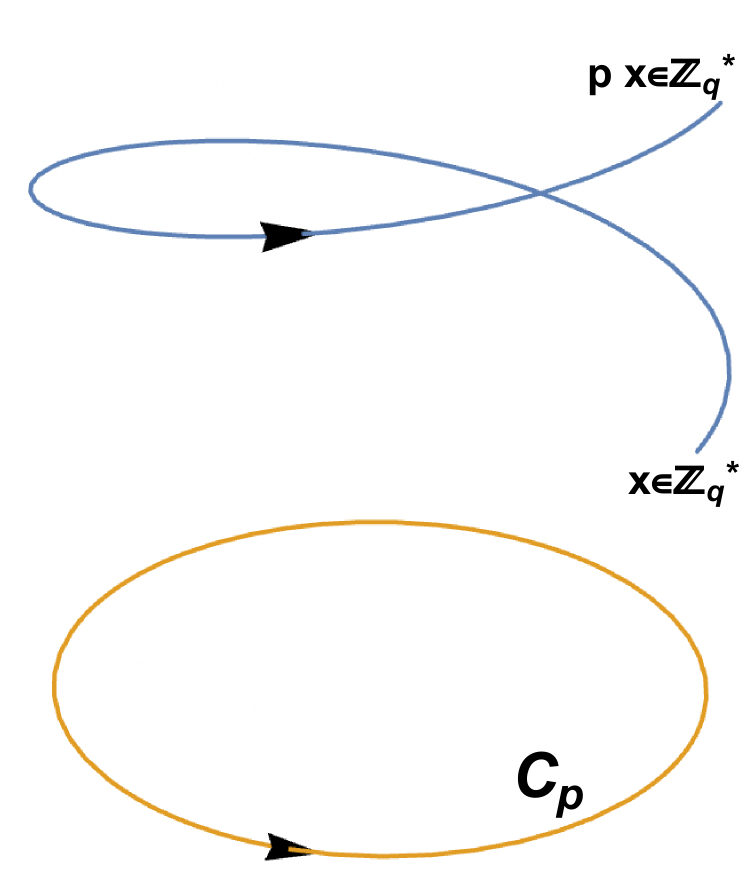}
\end{center}
\caption{Lifting the periodic orbit $C_p$\label{fig1} }
\end{figure}
Our second result gives the geometry underlying the analogy with the linking number of two knots $K,L\subset S^3$ in the three sphere. By definition, the linking number $\operatorname{lk}\left(K, L\right)$ is the monodromy  obtained by lifting the first knot $K$ in the maximal abelian cover $\left(S^3-L\right)^{ab}$ of the complement of the second knot $L$. Equivalently,  $\operatorname{lk}\left(K, L\right)$ is   the image in $\pi_1\left(S^3-L\right)^{a b}\simeq \mathbb{Z}$ of the canonical generator of $\pi_1(K)$
\begin{equation}\label{linking}
\ \ \pi_1(K) \longrightarrow \pi_1\left(S^3-L\right)^{a b}\simeq \mathbb{Z}.
\end{equation}
 In the analogy between knots and primes the role of the knots $K,L$ is played by two distinct primes $p,q$, while the  sphere is replaced by  $\Spec\, \Z$  and  $X_L:=S^3-L$ is replaced by  $\Spec\, \Z[1/q]$.\newline
Given two primes $p\neq q$ we consider the inverse image of the periodic orbit $C_p$ in the semilocal adele class space\footnote{the other places do not play a role}  associated to the set of  places $S=\{ p, q,\infty \}$. The semilocal adele class space is $\Gamma\backslash\left(\Q_p\times \Q_q\times \R\right)$ where $\Gamma:=\{\pm p^mq^n\mid m,n\in \Z\}$. 

\begin{thm}\label{mainintro} The  inverse image  of the periodic orbit $C_p$ in the semilocal adele class space $\Gamma\backslash\left(\Q_p\times \Q_q\times \R\right)$ associated to  $S=\{ p, q,\infty \}$ is the mapping torus of the canonical generator $\left\{F r o b_{p}\right\}$ acting by multiplication in the abelianized \'etale fundamental group  $\pi_1^{e t}(\Spec \, \Z[1/q])^{ab}$ of the complement of $q \in \Spec\, \Z$.	
\end{thm}
Thus  the periodic orbit $C_p$ plays the role of the first knot, and the monodromy of its lift  in the maximal abelian cover of the complement of $q \in \Spec\, \Z$ is the element $p\in \Z_q^*$. In the analogy between knots and primes developed in \cite{Mazur,Morishita,Morishita1}, $p\in \Z_q^*$  plays the role of the linking number\footnote{and allows one to compare the quadratic reciprocity with the antisymmetry of the linking number} of $p$ with $q$.\newline
Section \ref{sect3} provides a detailed description of the natural stratification of the semilocal adele class space $\A_{\Q,S}=\Gamma\backslash\left(\Q_p\times \Q_q\times \R\right)$ that generates a spectral sequence for the $K$-theory of the associated $C^*$-algebras.\newline
Finally, section \ref{sect4} describes the canonical 
 codimension one foliation of the three dimensional manifold $\Gamma\backslash\left(\A_{\Q,S}\times \underline{E \Gamma}\right) =\Gamma\backslash\left(\Q_p\times \Q_q\times \R \times \R^2\right)$   source of the Baum-Connes map \cite{BC} which computes   the $K$-theory of the $C^*$-algebra $\cA=C_0(\A_{\Q,S})\rtimes \Gamma$.

\section{Proof of Theorem \ref{mainglobalintro}}

Let $ \mathbb{A}^f$ be the ring of finite adeles of $\Q$. Let $(a,\lambda)\in \mathbb{A}^f\times \R= \A_\Q$ be an adele with $\lambda>0$. 
 Let $\Phi$ be the map from $\left(\mathbb{A}^f / \hat{\mathbb{Z}}^{\times}\right) \times \mathbb{R}_{+}^{\times}$ to subgroups of $\mathbb{R}$ defined by
$$
\Phi(a, \lambda):=\lambda H_a,\ \  H_a:=\{q \in \mathbb{Q} \mid a q \in \hat{\mathbb{Z}}\}.
$$
Then (\cite{CCas}, Lemma 3.7) $\Phi$ is a bijection between the subset of $X_\Q=\Q^\times \backslash \A_\Q/\hatz$  formed  of adele classes with non-zero archimedean component,   and the set of non-zero subgroups of $\mathbb{R}$ whose elements are pairwise commensurable. 
Let $\tilde \pi:\A_\Q\to\Q^\times \backslash \A_\Q\to X_\Q$ be the projection.
\[
\begin{tikzcd}
\A_\Q \arrow[rr]\arrow[dr,"\tilde\pi"'] & & \Q^\times \backslash \A_\Q\arrow[dl,"\pi"]\\
& X_\Q &
\end{tikzcd}
\]
Then  $\tilde\pi^{-1}(C_p)\subset \A_\Q$ is the saturation  $\Q^\times F$,  of the subset 
$$
F=\{a=(a_v), \ \ a_\infty >0, \ a_p=0, \ a_v\in \Z_v^* \qqq v\notin \{p,\infty\}\}\subset \A_\Q.
$$
The residual diagonal action of $\Q^\times$ on $F$ is given by  $p^\Z\subset \Q^\times$. Thus one obtains
\begin{equation}\label{inverseim}
\pi^{-1}(C_p)=\left(\R_+^*\times \prod_{q\neq p} \Z_q^*\right) /p^\Z
\end{equation}
where $p$ is diagonally embedded in the product. This quotient is by construction the mapping torus of the multiplication by $p$ in the compact group $\prod_{q\neq p} \Z_q^*$. We now interpret these terms in the language of \'etale abelianized fundamental groups. This derives from the following facts:

\begin{itemize}
\item The abelianized \'etale fundamental group\footnote{The abelianization takes care of the need to specify the base point} $\pi_1^{e t}(\Spec \, \Z_{(p)})$  is canonically isomorphic to $\prod_{q\neq p} \Z_q^*$.
\item The image $\pi_1^{e t}(r^*)\left\{F r o b_{p}\right\}$ in the étale abelianized fundamental group $\pi_1^{e t}(\Spec \, \Z_{(p)})\simeq \prod_{q\neq p} \Z_q^*$ is equal to $p$ diagonally embedded in $\prod_{q\neq p} \Z_q^*$ (see \eqref{residue}).
\end{itemize}
The first fact follows from \cite{Lenstra} (Corollary 6.17), together with the determination of the maximal abelian extension of $\Q$ in which $p$ is unramified as obtained by adjoining all roots of unity of order prime to $p$\footnote{following  the local to global proof of the Kronecker-Weber theorem}. Its Galois group is $\prod_{q\neq p} \Z_q^*$. The second fact follows since the action of $\left\{F r o b_{p}\right\}$ on roots of unity is given by raising to the power $p$.\vspace{.1in}

The above proof elucidates the geometric aspects of the well-known analogy between primes and knots through  the role of the scaling site $X_\Q$, its periodic orbits and their liftings to the adele class space of the rationals viewed as an abelian cover of $X_\Q$.

\section{The semilocal space $\Gamma\backslash \left(\Q_p\times \Q_q\times \R\right)$}\label{sect3}

\begin{figure}[H]
    \centering\vspace*{.3in}
    
\subfloat[]{{\includegraphics[width=4cm,valign=m]{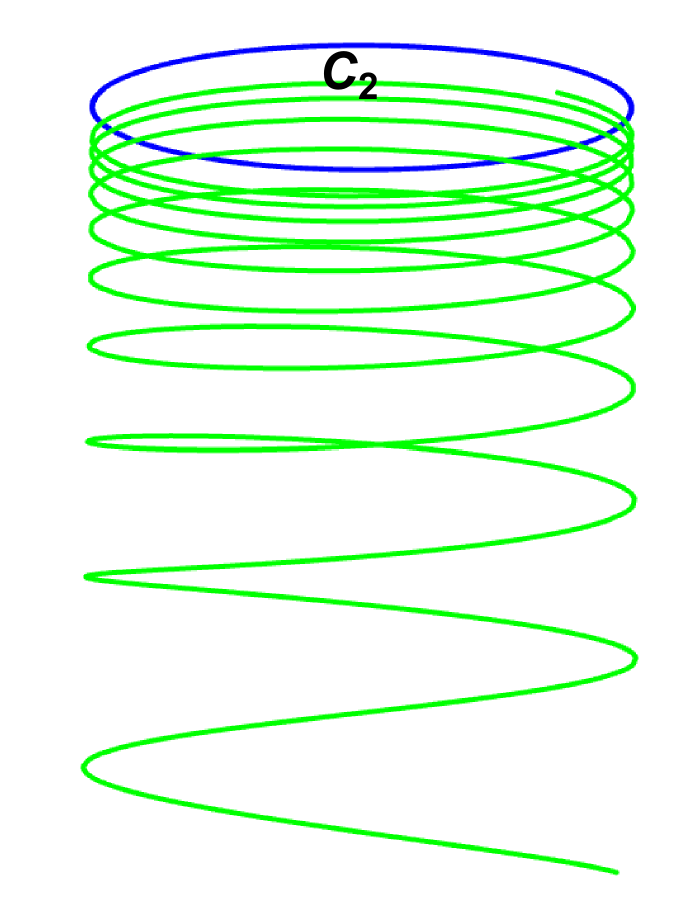} }}%
    \hspace{1in}
\subfloat[]{{\includegraphics[width=5.5cm,,valign=m]{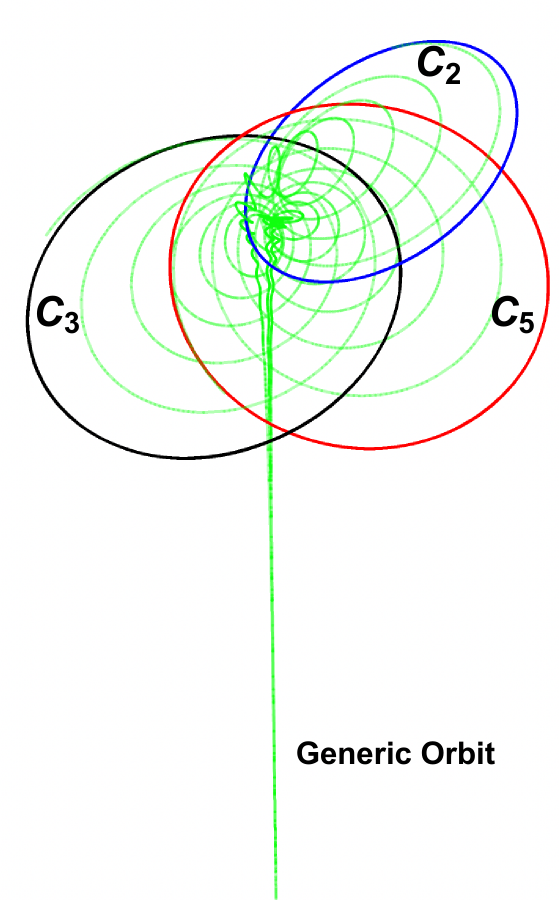} }}%
    \caption{(A) shows how the generic orbit (in green) for the action of $\R_+^*$ on $X_\Q$ is dense in the periodic orbit $C_p$. This holds for any $p$. (B) shows how the density of the generic orbit in every $C_p$ imitates the density of the generic point of $\Spec\,\Z$ (compare  with \cite{Manin}, page 16).}%
    \label{fig2}%
\end{figure}

Let $S$ be a finite set of places of $\Q$, with $\infty \in S$. We briefly recall the construction of the semilocal adele class space associated to $S$. One replaces the ring of adeles by 
the locally compact ring
$
 \A_{\Q,S}=\prod_{v\in S} \Q_v
$
in which  $\Q$  embeds diagonally. Let
$\cO_{\Q,S}$ be the subring of $\Q$ given by rational numbers whose
denominator only involves primes $p \in S$. In other words,
$
 \cO_{\Q,S}=\{q\in \Q\,|\, |q|_v\leq 1\,,\ \forall v\notin S\}\,.
$
The group $\Gamma$ of invertible elements of  $\cO_{\Q,S}$ is 
\begin{equation}\label{GL1QS}
\Gamma=\GL_1(\cO_{\Q,S})= \{ \pm p_1^{n_1} \cdots p_k^{n_k} \, :\,  p_j
\in S \setminus\{ \infty \} \,,\, n_j\in \Z\}.
\end{equation}
The semilocal adele  class space 
 is the quotient $X_{\Q,S}^{ab}:=\Gamma\backslash\A_{\Q,S}$. The map  
\begin{equation}\label{XQS}
\pi_S:X_{\Q,S}^{ab}\to X_{\Q,S}=\Gamma\backslash\A_{\Q,S}/\left(\prod_{v\in S\setminus \infty}\Z_v^*\right)
\end{equation}
 plays the role of the maximal abelian cover of the  double quotient $X_{\Q,S}$, the semilocal analogue  of $X_\Q$. 
The group 
\begin{equation}\label{JS} J_S=\Gamma\backslash\GL_1(\A_{\Q,S})\simeq \R_+^*\times  \prod_{v\in S\setminus \infty}\Z_v^*
\end{equation}
is the analogue of the idele class group and acts naturally on  $X_{\Q,S}^{ab}$ by multiplication.\newline
We let $p\neq q$ be two primes and  $S=\{ p, q,\infty \}$. The ring $\cO_{\Q,S}$
is $\Z[1/p,1/q]$ and the  abelianized \'etale fundamental group of its spectrum  is 
\begin{equation}\label{pi1}
\pi_1^{e t}(\Spec \, \cO_{\Q,S})^{ab}=\Z_p^*\times \Z_q^*.
\end{equation}
Next we describe the natural stratification of the space $X_{\Q,S}^{ab}$.
\begin{prop}\label{ergodic23} Let $p\neq q$ be two primes,  $S=\{ p, q,\infty \}$ and $\pi_S:X_{\Q,S}^{ab}\to X_{\Q,S}$ as in \eqref{XQS}. Then \newline
$(i)$~The orbits of the action of the group $J_S$ on $X_{\Q,S}^{ab}$ are indexed by the subsets of $S=\{ p, q,\infty \}$ as follows 
$$
\Omega_Z:=\Gamma\backslash\{(a_v)\in \A_{\Q,S}\mid a_v= 0\qqq v \in Z, \ a_v\neq 0\qqq v \notin Z\}\qqq Z\subset S.
$$
$(ii)$~The orbit  $\Omega_{\{p\}}$ associated to the subset $\{p\}\subset S$ is  the inverse image $\pi_S^{-1}(C_p)$ of the periodic orbit  $C_p$ of length $\log p$ (and a similar result holds for the subset $\{q\}\subset S$).\newline
$(iii)$~The orbit $\pi_S^{-1}(C_p)$ is equivariantly isomorphic to the mapping torus of the multiplication by $p$ in the compact group $ \Z_q^*$.	
\end{prop}
\proof $(i)$~Follows from the transitivity of the action of $\GL_1(\Q_v)$ on $\Q_v\setminus \{0\}$ for any  $v$.\newline
$(ii)$~The periodic orbit $C_p$ of length $\log p$  for the action  of the scaling group $\R_+^*$ in $X_{\Q,S}$ is equal to $\pi_S(\Omega_{\{p\}})\subset X_{\Q,S}$ and is identical to the quotient $p^\Z\backslash\R_+^*$. This description shows that $C_p$ appears in the boundary of the free orbit $\pi_S(\Omega_{\emptyset})\subset X_{\Q,S}$ as shown in Figure \ref{fig2}.  \newline
$(iii)$~ One has $\Omega_{\{p\}}=\pi_S^{-1}(C_p)$. The action of $\{\pm q^n\}$ on $\{0\}\times \left(\Q_q\setminus \{0\}\right)\times \R^*$ admits $\{0\}\times \Z_q^*\times \R_+^*$ as fundamental domain. This domain is preserved under the action of $p^\Z$. The quotient 
$$
p^\Z\backslash\left(\Z_q^*\times \R_+^*\right)
$$
is the mapping torus of the multiplication by $p$ in the compact group $ \Z_q^*$.
\endproof
There are $8$ orbits $\Omega_Z$,  $Z\subset S=\{p,q,\infty\}$.  Proposition \ref{ergodic23} describes two of them, the remaining six orbits are described as follows:
\begin{enumerate}
     \item The orbit $\Omega_{\emptyset}\subset X_{\Q,S}^{ab}$  is the idele class group $J_S$ of  \eqref{JS}.
   \item The orbit $\Omega_{\{\infty\}}\subset X_{\Q,S}^{ab}$ is the quotient $\left(\Z_p^*\times \Z_q^*\right)/\pm 1$. 
   \item  The orbit $\Omega_{\{p,q\}}$ is the quotient $(p^\Z q^\Z)\backslash\R_+^*$   \item  The orbit $\Omega_{\{p,\infty\}}$ is the quotient  of $ \Z_q^*/\{\pm 1\}$ by the action of $\Z$ by multiplication by powers of $p$. 
   \item  The orbit $\Omega_{\{q,\infty\}}$ is the quotient  of $ \Z_p^*/\{\pm 1\}$ by the action of $\Z$ by multiplication by powers of $q$. 
   \item The orbit  $\Omega_{\{p,q,\infty\}}$ is a single point.
\end{enumerate}
To be more precise one can describe for each orbit the corresponding $C^*$-algebra up to Morita equivalence. For $(3)$ one obtains the noncommutative torus associated to the irrational number $\log p /\log q$. For $(4)$ one gets the cross product  $
C(\Z_q^*/\pm 1)\ltimes p^\Z
$
and similarly for $(5)$. Finally for $(6)$ one gets the group ring $C^*(\Gamma)$. By construction the group $J_S$ acts by automorphisms on these $C^*$-algebras.
\vspace{-.2in}

\subsection*{Proof of Theorem \ref{mainintro}}
This follows from Proposition \ref{ergodic23} $(iii)$, together with the following facts :
\begin{itemize}
\item The abelianized \'etale fundamental group $\pi_1^{e t}(\Spec \, \Z[1/q])^{ab}$  is canonically isomorphic to $\Z_q^*$.
\item The image $\pi_1^{e t}(r^*)\left\{F r o b_{p}\right\}$ is equal to $p\in  \Z_q^*$.
\end{itemize}
As in the proof of Theorem \ref{mainglobalintro}, the first fact follows using \cite{Lenstra} (Corollary 6.17), together with the determination of the maximal abelian extension of $\Q$ only ramified at $q$ (obtained by adjoining all roots of unity of order a power of $q$). Its Galois group is $\Z_q^*$. The second fact follows since the action of $\left\{F r o b_{p}\right\}$ on roots of unity is given by raising to the power $p$.

\section{The classifying space $\Gamma\backslash\left(\A_{\Q,S}\times\underline{E \Gamma}\right)$ and its codimension $1$ foliation} \label{sect4}

  Let $S=\{ p, q,\infty \}$ as above. The $C^*$-algebra $\cA=C_0(\A_{\Q,S})\rtimes \Gamma$ of the noncommutative space $X_{\Q,S}^{ab}:=\Gamma\backslash\A_{\Q,S}$ is not of type I (see \cite{CMbook}, Lemma 2.28 Chapter II).  The  understanding of $\cA$ at the level of $K$-theory is obtained by considering the assembly map \cite{BC} computing the $K$-theory of the  $C^*$-algebra $\cA$ in terms of  the universal  proper action  $\underline{E \Gamma}$ of $\Gamma$. Since  $\Gamma$ is commutative the various definitions of the cross product $C^*$-algebra all coincide with the reduced cross product, thus we drop the lower index $r$ in the notation $C^*_r$. Moreover the assembly map is an isomorphism. One obtains in this way  an ordinary space $Y=Y_{\Q,S}=\Gamma\backslash\left(\A_{\Q,S}\times\underline{E \Gamma}\right)$ which resolves the noncommutative nature of the quotient $\Gamma\backslash\A_{\Q,S}$. The group $\Gamma$ is $\Z^2\times\{ \pm 1\}$ and the universal  space for proper actions $\underline{E \Gamma}$ is $\R^2$ where $\Z^2$ acts by translations and $\{ \pm 1\}$ acts trivially. Thus one obtains  the following three dimensional space, quotient by the diagonal action of $\Gamma$,
 \begin{equation}\label{3space}
 Y=\Gamma\backslash\left(\A_{\Q,S}\times \underline{E \Gamma}\right) =\Gamma\backslash\left(\Q_p\times \Q_q\times \R \times \R^2\right).
 \end{equation}
 \vspace{.1in}
 
 We stress that the dimension three in the next proposition   is due to taking a pair of primes and has no relation to the homological dimension of $\Spec\,\Z$.
  \begin{prop}  	
 \label{bciso} $(i)$~The space $Y$ is locally compact of dimension $3$.\newline
  $(ii)$~The action of $\R^2$ by translations  on the   space $Y$ defines a codimension $1$  lamination of $Y$ whose space of leaves is $X_{\Q,S}^{ab}$.
    \end{prop}
    \proof  $(i)$~The properness of the action of $\Gamma$ on $\underline{E \Gamma}$ implies  that its action on $\left(\A_{\Q,S}\times \underline{E \Gamma}\right) $ is also proper. Thus $Y$ is Hausdorff  and since $\left(\A_{\Q,S}\times \underline{E \Gamma}\right) $ is locally compact, $Y$ is  locally compact. The dimension of $\Q_p\times \Q_q\times \R \times \R^2$ is equal to $3$ and is unaltered when passing to the quotient.\newline
    $(ii)$~By construction the orbits of the action by translation of $\R^2$ on $\left(\Q_p\times \Q_q\times \R \times \R^2\right)$ define a lamination which is invariant under the action of $\Gamma$ and hence descends to the proper quotient $Y$. The space of leaves of the first lamination is $\left(\Q_p\times \Q_q\times \R \right)$ and its quotient by $\Gamma$ is $X_{\Q,S}^{ab}$.\endproof 
     The simplest meaningful computation of the $K$-theory of the involved $C^*$-algebras  is for the cross product $A$ associated to the union in $X_{\Q,S}$ of the generic orbit with the three periodic orbits $C_p,C_q,C_\infty$. One obtains that $K_0(A)\simeq \Z^3$ reflects the presence of the three periodic orbits, while $K_1(A)\simeq \Z^2$ reflects the one-dimensionality of the periodic orbits $C_p,C_q$.
\begin{rem}    
$(i)$~The theory of the \'etale fundamental group for schemes can be viewed as the extension of Galois theory to schemes (\cite{Lenstra}). The classical abelian class field theory provides an isomorphism of the Galois group of the maximal abelian extension of a global field with an adelic group. It is natural to wonder if this 
 isomorphism can be extended to compare the behavior of the \'etale fundamental group of schemes with an adelic geometric construction.
Theorems \ref{mainglobalintro} and \ref{mainintro} 
answer positively to this question for  schemes related to $\Spec\,\Z$ in terms of the geometry of the adele class space viewed as a covering $\pi:X_\Q^{ab}\to X_\Q$ of the scaling site. This extension of the class field theory isomorphism to schemes ought to play a role in the understanding of  the subtle relation between the scaling site and $\Spec\,\Z$.	\newline

$(ii)$~The scaling site which plays a central role in the present paper admits a natural complexification $\mathscr C_\Q$ \cite{CCsurvey}, whose role is to pass from characteristic one to characteristic zero. It replaces the real half-line by the pro-\'etale cover of the punctured unit disk in the complex plane. Then $\mathscr C_\Q$  is interpreted as the moduli space of elliptic curves with a triangular structure. The extrapolation of the findings of the present paper to the context of 
$\mathscr C_\Q$
  has the potential to formulate the results in   crystalline terms in place of the \'etale one. 
\end{rem}


\end{document}